# Pareto-optimal treatment of uncertainties in model-based process design and operation

Dr. Jan Schwientek[1], Dr. Katrin Teichert[1], Jan Schröder[1], Dr. Johannes Höller[1], Dr. Norbert Asprion[2], Dr. Pascal Schäfer[2], Dr. Martin Wlotzka[2], and Prof. Dr. Michael Bortz[1,*]

Model-based process design and operation involves here-and-now and wait-and-see decisions. Here-and-now decisions include design variables like the size of heat exchangers or the height of distillation columns, whereas wait-and-see decisions are directed towards operational variables like reflux and split ratios. In this contribution, we describe how to deal with these different types of decisions in a multicriteria framework, offering an adjustability for the wait-and-see variables while at the same time respecting optimality guarantees on process KPIs.

**Author affiliations**

[1]Fraunhofer ITWM, Fraunhofer Platz 1, 67663 Kaiserslautern, Germany.

[2]BASF SE, Carl-Bosch-Strasse 38, 67056 Ludwigshafen am Rhein, Germany.

Email corresponding author: michael.bortz@itwm.fraunhofer.de

ORCID iDs of the authors

https://orcid.org/0000-0001-9740-3233 (Jan Schwientek)

https://orcid.org/0000-0002-2293-407X (Katrin Teichert)

https://orcid.org/0009-0003-8337-0807 (Jan Schröder)

https://orcid.org/0000-0003-3431-7818 (Johannes Höller)

https://orcid.org/0000-0002-7728-0499 (Norbert Asprion)

https://orcid.org/0000-0002-3268-8976 (Pascal Schäfer)

https://orcid.org/0000-0003-1794-4799 (Martin Wlotzka)

https://orcid.org/0000-0001-8169-2907 (Michael Bortz)

## 1    Introduction

Process design is a multicriteria, two-stage optimization problem. In the first stage, design-related variables, henceforth called here-and-now-variables (HNV), are fixed. These determine the kind and the size of unit operations. In the second stage, operational variables, called wait-and-see-variables (WSV), are chosen. Among those are, for example, temperatures, reflux and split ratios. The problem is multicriteria because the choice of variables is always governed by the effort to balance competing objectives [1]: For example, CAPEX and OPEX of a production process shall be minimized, while at the same time keeping high quality and security standards. Another example is minimizing the total annualized costs and the environmental impact.

In practice, often neither input conditions like the amount and concentrations of feed streams or the prices for energy and raw materials, nor process requirements like conditions on product purities are fixed to single values. Rather, a range of different conditions and requirements has to be met. Conceptually, the same applies for physical model parameters: These are known only within a certain confidence interval. The different values that may occur for the input conditions, process requirements and physical parameters are called scenarios in the following.

Here, this situation is addressed by seeking for values of HNV and WSV so that irrespective of the scenarios within their ranges, one can always guarantee that a particular choice of HNV and WSV realizes a best compromise between the conflicting objectives. Additionally, feasibility shall be guaranteed by fulfillment of all constraints in any scenario. This strategy makes sure that there is no hidden potential for further improvement in the objectives, unless the above-mentioned ranges or other restrictions are changed.

This work focusses on how to obtain suitable values for HNV and WSV in the multicriteria context of chemical process design: When deciding for a HNV, it is assumed that one knows the possible scenarios that may happen. Thus, HNV must be chosen in a way that all these possible scenarios are accounted for before one particular scenario realizes. This is known as robust optimization [2], [3]. However, when the process is being operated, it is assumed that it is possible to react to the realized scenario. This adjustability of decisions involving WSV is well known for single-criteria optimization problems [4] and leads to optimization potential and, thus, an improvement of the objectives with respect to the more conservative approach applied to the HNV. In this contribution, we will describe how to realize this adjustability in a multicriteria setting and apply it in the context of chemical process design and operation.

From the process engineer's perspective, this is a practical procedure to find designs and operating conditions to optimally hedge against typically occurring scenarios. Technically, this leads to multicriteria adjustable robust optimization problems, which belong to bilevel or even multilevel optimization.

The need to deal with uncertainties in the design of chemical production processes has been made clear in [5]. In this work, a stochastic programming approach has been chosen; details on the mathematical background can be found, e.g., in [6]. Its application to model-based chemical process design also includes the development of dedicated software modules [7]. Within the robust approach, quite some attention has been paid to the flexibility index, which is a measure of how large the ranges of the different scenarios can be [8],[9].

As to multicriteria optimization under uncertainties, Plazoglu and Arkun [10] were among the first to combine robustness considerations with a multicriteria optimization (MO) formulation. They introduced a two-stage approach in which the MO problem is decomposed into a series of structurally fixed single-objective problems using scalarization techniques—specifically, the ϵ-constraint method—while addressing uncertainty separately in each optimization. Similarly, Dantus and High [11] tackled the trade-off between cost and environmental impact using a stochastic optimization approach (simulated annealing), where Monte-Carlo simulations were employed to account for uncertainties across multiple process evaluations.

A generalized version of the two-stage approach was later proposed by Fu et al. [12] wherein a master MO algorithm defines a sequence of single-objective problems that incorporate uncertainty through sampling techniques, such as Hammersley sequence sampling. A similar framework, coupled with specialized genetic algorithms, was later presented by Kheawhom et al. [13], [14]. Kim and Diwekar [15] extended this approach to solvent design problems, where discrete optimization variables define solvent structures. Their study highlighted the impact of the probability distribution of uncertain parameters on optimization outcomes.

Beyond the two-stage approach, Chakraborty et al. [16], [17] applied MO under uncertainty to plant-wide optimization using superstructures. They introduced a flexibility index to assess design robustness, using it to evaluate the Pareto set obtained from a MO problem balancing cost against global warming potential. Similarly, Hoffmann et al. [18] initially computed a Pareto set under deterministic assumptions and subsequently assessed the effect of parameter uncertainties by sampling their distributions and visualizing objective deviations.

Another approach, parametric MO, was introduced by Hugo and Pistikopoulos [19] wherein multiple discrete scenarios are defined, and the MO problem is solved for each. The resulting Pareto sets are then compared. Tock and Marechal [20] demonstrated a probabilistic method using Monte-Carlo simulations to quantify the probability distribution of each Pareto set, further refining the uncertainty analysis within MO frameworks.

For examples of other kinds of robustness in a multicriteria setting the interested reader is referred to [21].

The remainder of this article is organized as follows: In section 2, the concepts of robustness and adjustability to hedge against uncertainties by choosing HNV and WSV in a multicriteria context are presented. Section 3 contains the application of these concepts to a toy example and to a binary distillation column. The article ends with an outlook.

## 2 Here-and-now vs. Wait-and-see variables in multicriteria optimization

Addressing process design tasks as multicriteria optimization (MO) problems has proven highly effective. For a comprehensive overview, readers are directed to references [1], [22],[23]. Integrating MO within an interactive decision-support framework enables the exploration of the Pareto boundary, uncovering Pareto solutions that single-objective methods miss [1]. A MO problem is written as

$$\min_{x \in X, y \in Y} F(x, y, u_{nom}) \qquad (1)$$

where $x$ ($y$) denotes the HNV (WSV) with their feasible domains $X$ ($Y$), respectively. These feasible domains may be given explicitly or implicitly by respecting additional inequality or equality constraints, not given explicitly above. The vector of functions $F(x, y, u_{nom}) = \left(f_1(x, y, u_{nom}), \ldots, f_M(x, y, u_{nom})\right)^T$ contains $M$ objective functions. The model parameters are fixed to their nominal values $u_{nom}$. The solution of the MO problem (1) is the nominal Pareto set, which is of high value for transparent decision support in process design and optimization [1], [24].

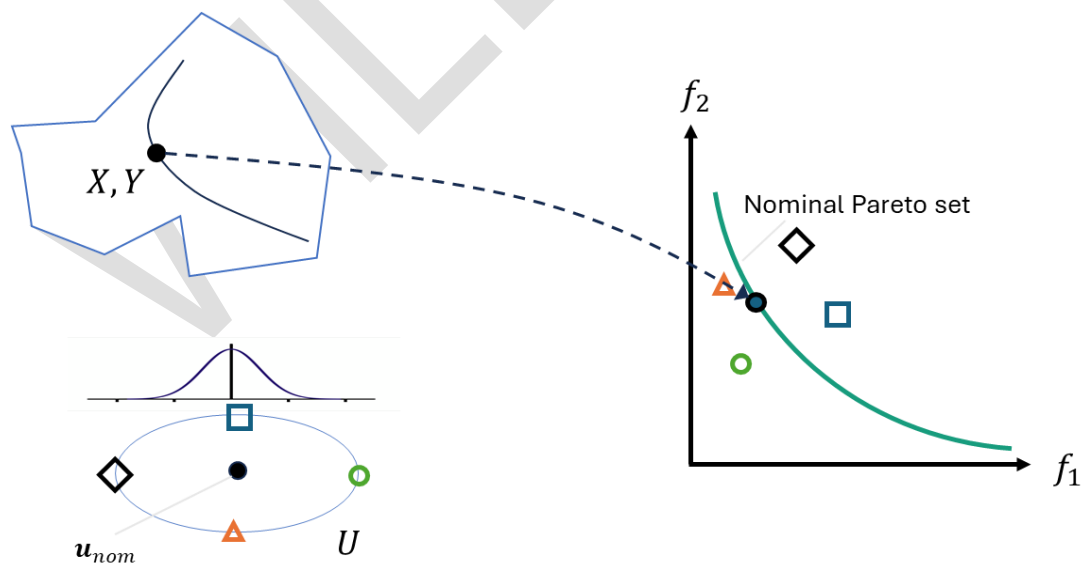

**Figure 1.** Illustration of the nominal Pareto set for a bi-objective optimization problem. Scenarios have been chosen according to a one-at-a-time strategy in parameter space. Here, this has been motivated by a hypothetical elliptic confidential region, which could result from parameters normally distributed around their nominal values.

In order to investigate the impact of parametric uncertainties on the nominal Pareto set, it is helpful to do a sensitivity analysis [25] for those parameters that can differ from their nominal values $u_{nom}$. Therefore, one starts with an assumption about the relevant set $U$ of scenarios deemed to occur with a realistic probability. Then, for a given Pareto point, one calculates the values of the objective functions, as shown in Figure 1. Using a scattering measure, one may measure the impact of different scenario values on the Pareto set [25]. Typically, certain scenarios lead to worse values in the objective functions, compared to the nominal Pareto set (in Figure 1, these are the orange and black scenarios).

Among those, one or more scenarios can be identified to lead to the worst-case in objective space. It should be noted that this "worst-case" refers to the chosen scenarios in $U$. Here it is assumed that these scenarios are realistic alternatives. It goes without saying that the larger the set $U$, the more conservative the component-wise robust Pareto set will be, and the higher the cost compared to the nominal Pareto set.

If one thinks of the Pareto set as benchmark showing the best compromises attainable for all scenarios, the original MO problem (1) must be modified. This modification is done in two steps.

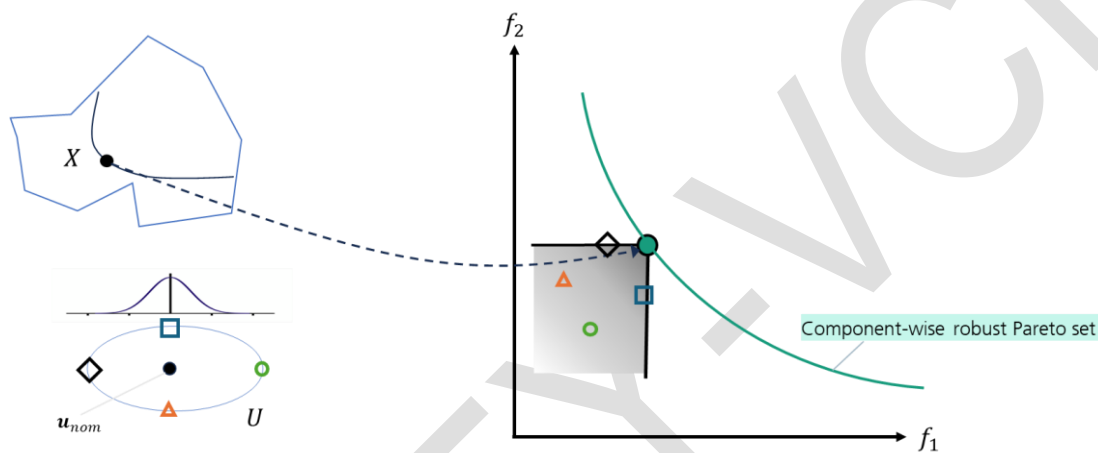

**Figure 2.** Illustration of the component-wise robust Pareto set compared to the nominal case with the same scenarios as in Figure 1.

## 2.1 Component-wise robust Pareto set

In a first step, consider the robust counterpart of the nominal MO problem (1) [2],

$$\min_{x \in X, y \in Y} \max_{u \in U} F(x, y, u), \qquad (2)$$

which is called robust MO (RMO) formulation in the ongoing. The inner optimization problem means that for each objective function, i.e., each component of $F$, one is looking for the worst case that may happen under all $u \in U$. These may differ for the different objective functions [26]. One then determines one value for the HNV $x$ so that these worst cases are minimized. The formulation (2) leads to a component-wise robust Pareto set. It is of bilevel structure, tractable by techniques of semi-infinite programming [27] like, for example, by the algorithm of Blankenship and Falk [28]. Since we have in mind to apply the approach to flowsheet simulations with large equation systems, the inner optimization problem is replaced by the maximization over the finite set of discrete scenarios in $U$, introduced above. Thus, the solution of the inner problem amounts to an enumeration over all scenarios. The resulting component-wise robust Pareto set is illustrated in Figure 2.

The RMO (2) is perfectly suited to deal with the HNV $x \in X$: At the stage when those are fixed, one has a certain idea about what realistically may happen (the set $U$) and wants to hedge against these scenarios. This also means that in the multicriteria context, the resulting Pareto points are in so far pessimistic as they do not need to correspond to the actually occurring worst cases, as shown in Figure 2: Compared to the full green circle solution, the blue solution is better in $f_2$, the black solution is better in $f_1$. However, there is no worst-case solution, for a given $x$, that is simultaneously better in both objectives, as expected.

As far as the WSV $y \in Y$ are concerned, the component-wise RMO problem (2) is too pessimistic: Having to decide on a WSV means that the scenario at hand is known, so that one can directly react to it. How this reaction looks like, that is, how the value of the WSV is chosen, can differ between scenarios. This means that in contrast to handling HNV, one does not seek values for $y$ that are the same for all scenarios: The HNV are adjustable to the scenarios.

## 2.2 Adjustable robust Pareto set

To exploit this adjustability, one proceeds as follows: For each scenario $u \in U$, the Pareto set with respect to the WSV $y$ for fixed values of the HNV $x$ is determined. From this, the adjustable robust Pareto set is extracted. The values for $x$ are chosen such that with respect to them, component-wise robustness is respected. This results in the following formulation:

$$\min_{x \in X} \max_{u \in U} \min_{y \in Y} F(x, y, u) \tag{3}$$

which is denoted as multicriteria adjustable robust optimization problem (MARO) in the following.

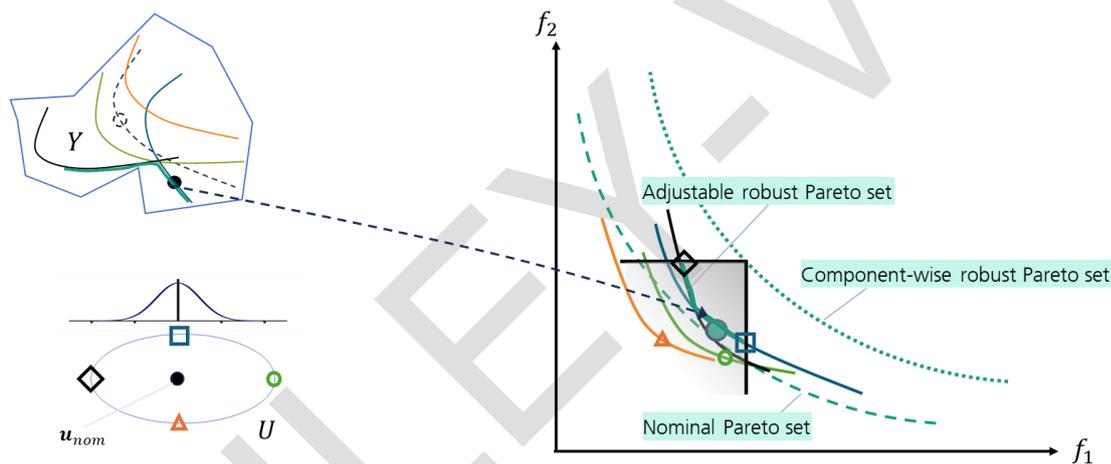

**Figure 3.** Illustration of the adjustable robust Pareto set compared to both the nominal and the component-wise robust Pareto sets with the same scenarios as in Figure 1.

Again, for a practical realization, this is done for finitely many scenarios, as illustrated in Figure 3. It should be noted that the adjustable robust Pareto set is expected to dominate the component-wise robust set, reflecting the effect of adjustability of the WSV. As can be seen from Figure 3, the Pareto set for each scenario is calculated individually, and then the resulting worst-case set is identified. The HNV are variables are chosen so that all possible scenarios for the resulting values of the objective functions are covered.

The conclusion to be drawn is that accounting for the adjustability of the WSV allows to unveil additional optimization potential that is not accessible if both WSV and HNV are treated on an equal footing, oriented towards component-wise robustness.

The formulation (3) is a challenging trilevel optimization problem. To be able to solve it with a flowsheet simulator, i.e., as usual nonlinear optimization problem, it is assumed that the maximization over $u \in U$ is done for a finite set of fixed, pre-defined scenarios $u_{i=1,\dots,K}$, as sketched in Figures 1-3. To take account of the adjustability of the WSV $y \in Y$, a replication is done: For each $u_i$, an optimization variable $y_i$, representing the function value of the unknown solution function $y(u_i)$, is introduced. This replication and the discrete scenarios allow to reformulate (3) as follows:

$$\min_{\substack{x \in X, \\ y_{1,\ldots,K} \in Y^K}} \tilde{F}(x, y_{1,\ldots,K}), \tag{4}$$

with $\tilde{F}(x, y_{1,\ldots,K}) = \left(\max_{k=1,\ldots,K}\{f_1(x, y_k, u_k)\}, \ldots, \max_{k=1,\ldots,K}\{f_M(x, y_k, u_k)\}\right)^T$, which is a computational tractable variant of (MARO). This is a standard MO problem for the component-wise maximum of objectives $\tilde{F}$ meaning that the selection of the maximum has to be done in each step of the optimization algorithm. Thus, a standard nonlinear solver can be used to solve (4). This approach will be pursued in the following.

An alternative way to the above sketched replication is to use a parametric approximation of the unknown function $y(u)$ in (3), a so-called decision rule, leading to [29]

$$\min_{x \in X, y(u) \in Y} \max_{u \in U} F(x, y(u), u). \tag{5}$$

In order to address this problem, one may think of a Taylor expansion in $u$ around the nominal scenario, i.e., up to linear order, $y = a + B u$ with unknown parameter vector $a$ and matrix $B$. In such an approach, these are optimization variables; for details of this linear affine decision rule, cf [30], [31]. In this work, the replication formulation (4) of the MARO problem will be used.

Before proceeding to an application example, the following comments are in place: The first comment refers to the tradeoff between robustness and the original objectives. Since the robust Pareto sets are typically dominated by the nominal Pareto set it is natural to ask for the costs of robustness. Therefore, corresponding scenarios should be compared: Knowing the points in $X, Y$ leading to the robust Pareto sets, one can calculate the corresponding objective values for the nominal scenario. This is to be compared with the nominal Pareto set, thus yielding the costs of robustness in each objective function.

The second comment addresses the choice of scenarios. The larger the set of scenarios, the larger the computational effort. Above, scenarios have been chosen according to a presumed Gaussian probability distribution with an ellipsoidal joint confidence region. In practical application, the process engineer may have scenarios in mind that should be considered. Ideally, one knows which scenarios lead to a worsening of the objective values: These are the only ones that enter the above approach. For the component-wise robust Pareto set, it would be sufficient to know, for each $x$, the true worst-case scenario $u_x^*$ out of the continuous set $U$. It can be determined by applying techniques of semi-infinite programming to the problem (1) with a continuous set $U$. This goes beyond the scope of this paper and will be investigated in future work.

The third comment is directed towards the choice of $U$. Throughout this work it is assumed that at least certain typical scenarios in $U$ are known. Experience in practice shows though that there are situations where it is more realistic to assume that certain bounds on the objectives exist, and that the interest is in maximizing the set of scenarios that can be covered by a robust approach while respecting these bounds. This perspective is known as inverse robust optimization [32], solvable by an extension of the Blankenship-Falk algorithm [33], and offers a direction for future research as well.

## 3   Application examples

In this section, two application examples will be discussed: One simple toy example and one industrial application, namely a distillation column for the separation of a binary mixture. The results presented for the industrial example have been obtained by implementing the approach (4) in Chemasim, BASF's in-house flowsheet simulator. To solve the optimization problem, Schittkowski's suite with SQP solvers has been used [34].

**3.1 Toy example**

Consider two objective functions, both to be minimized,

$$f_1(t, u) = 1 - \cos(t) + u$$
$$f_2(t, u) = 1 - \sin(t) - u$$

with $u \in [-0.1, 0.1]$, $t \in [0, \frac{\pi}{2}]$. In this subsection, the difference of treating $t$ as a HNV (i.e. $t \equiv x$ in the notation of section 2), leading to an RMO problem, compared to treating it as a WSV within MARO (which means $t \equiv y$) will be highlighted. The results are shown in Fig. 4 below.

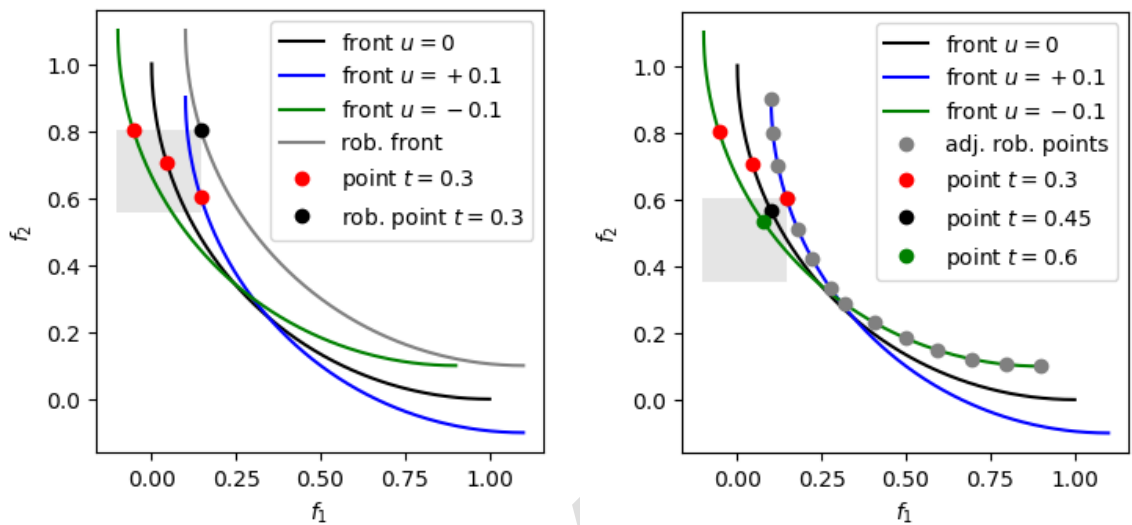

**Figure 4.** Robust (left hand side) and adjustable robust (right hand side) Pareto sets for the toy example.

In Fig. 4, on the left-hand side, the robust Pareto set is given, together with the identification of one of its points for a given value of the HNV $x$. On the right-hand side of Fig. 4, the MARO Pareto set is shown, showing the benefit of adjustability for the WSV $y$.

**3.2 Distillation column**

As mentioned above, the results in this subsection have been obtained by implementing the approach (4) in BASF's inhouse flowsheet simulator Chemasim.

A standard distillation column with a reboiler at the bottom and a condenser at the top with bottom and top production streams is considered to separate a binary zeotropic mixture of Methanol (MeOH) and Methyl formate (MF). The light boiler MF is withdrawn at the top. The challenge consists in identifying HNV and WSV in order find the Pareto set for minimizing CAPEX and OPEX, assuming uncertainties in the load of the column, the feed composition and the activity coefficients. The load of the column $l$ is defined as feed mass flow rate divided by 8000 kg/h; the uncertainty of the feed composition is captured by the mass fraction $w_{MF}$ of MF in the feed and the uncertainty in the activity coefficients is modeled by the approach described in [35]. Their ranges are:

$$l \in [0.6, 1.2] \text{ with nominal value } l^{nom} = 1.0$$
$$w_{MF} \in [0.78, 0.82] \text{ kg/kg with nominal value } w_{MF}^{nom} = 0.8 \text{ kg/kg}$$
$$F_{12} \in [0.9, 1.1] \text{ with nominal value } F_{12}^{nom} = 1.0,$$

where $F_{12}$ is a prefactor multiplying the limiting activity coefficient of MeOH in MF. The HNV considered here include the number of equilibrium stages $N$, the feed stage $N_f$, the diameter $D$ of the column and the areas of the heat exchangers in the reboiler $A_r$ and condenser $A_c$. Two WSV are considered, namely the reflux ratio $R_V$ and the specific heat duty $\dot{Q}_r$ of the reboiler. The objective functions, namely CAPEX and OPEX, are calculated by realistic cost functions.

|  | Nominal (MO) | Robust (RMO) | Adjustable robust (MARO) |
|---|---|---|---|
| $N \in [10, 150]$ | 13 - **150** | 19 - 80 | 13 - **150** |
| $N_f \in [3, 40]$ | 4 – 6 | 4 - 6 | 4 - 5 |
| $D \in [0.8, 2.0]$ in m | 0.95483 – 0.99768 | 1.1615 – 1.6909 | 1.0454 – 1.5038 |
| $A_r \in [50, 1000]$ in m$^2$ | 150.8 – 164.38 | 220.1 – 393.09 | 181.34 – 894.3 |
| $A_c \in [50, 1000]$ in m$^2$ | 132.72 – 145.39 | 196.48 – 758.16 | 159.62 – 693.99 |
| $R_V \in [0.5, 2.0]$ | 0.59099 – 0.74287 | 0.9762 – 1.2153 | 0.58709 – 0.89498 |
| $\dot{Q}_r \in [0.0625, 0.375]$ in kWh/kg | 0.1885 – 0.20547 | 0.22918 – 0.25545 | 0.1882 – 0.22184 |

**Table 1.** Ranges of optimization variables found in the optimization runs for the MF-MeOH distillation column. If bounds of these ranges coincide with lower/upper bounds set beforehand for the optimization, they are printed in bold. Variables above (below) the thick line are HNV (WSV).

Table 1 summarizes the ranges of the optimization variables (both HNV and WSV) for calculating the nominal, robust, and adjustable robust Pareto sets.

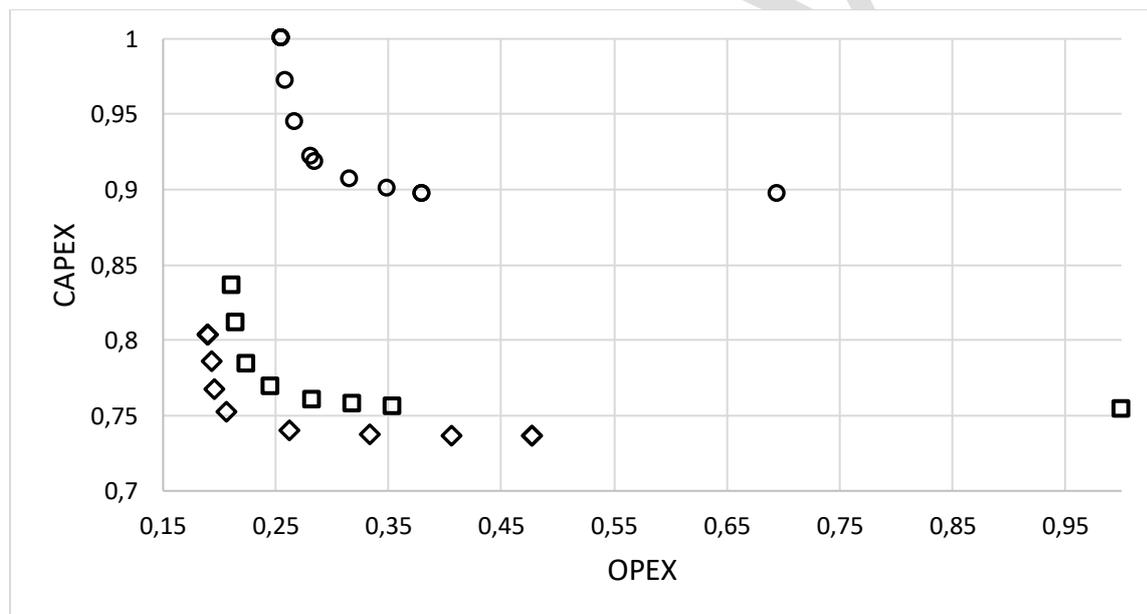

**Figure 5.** Pareto sets obtained for the distillation column separating MF and MeOH. Diamonds denote the nominal Pareto set; circles the robust and squares the adjustable robust Pareto sets. Units of the axes have been chosen such that the maximal value of each objective has been set to 1.

The result is shown in Figure 5. One can nicely observe the differences between the nominal, robust, and adjustable robust Pareto sets. Especially, it is worthwhile observing that by taking the adjustability of the WSV into account, hedging against the uncertainties is possible by only a small compromise compared to the nominal Pareto set. The robust and adjustable robust optimizations were carried out on the basis of one-at-a-time scenarios. As far as computation times are concerned, the nominal Pareto set requires ~ 5 s per point, the robust Pareto set ~ 2 min per point and the adjustable robust Pareto set ~ 15 min per point.

As illustrated above, the robust Pareto front is far more conservative. In this example, however, this is not due to the component-wise WC of the objective functions because only the operating costs depend on the uncertainties. Rather, this is due to the constraints which are contained in the problem:

- one upper bound each for the heating and cooling capacity of the reboiler resp. the condenser,
- one minimum purity requirement each for the products MeOH and MF, and
- a restriction on the F-factors of the distillation column (not to be confused with $F_{12}$ used above to model the uncertain activity coefficient)

For instance, the uncertain MF concentration in the feed leads to different scenarios (with high and with low concentration of MF) limiting the competing purity requirements for the separated substances, i.e., being the worst cases. This reduces the set of feasible decision variables and, thus, also leads to a more conservative Pareto set: In the RMO (circles), we require that with a chosen design, there is one reflux ratio and one specific duty, so that the separation stays feasible for any scenario. In particular, as we want to cope with both lowest MF content and highest with the same set of variables (HNV and WSV), we generally require much larger distillation columns. In contrast, in the MARO (squares), we can freely select the WSV dependent on the realized scenario. That is, here, we can adjust the reflux ratio and the specific duty to account for a changing MF content while still meeting the specifications and thereby mitigate the need for very high numbers of separation stages.

## 4     Conclusion

Different concepts for handling parametric uncertainties in flowsheet simulation und optimization have been developed since the seminal work [5]. Especially, significant progress has been made towards handle uncertainties in a multicriteria framework. This now makes it possible to incorporate the concepts suitable for industrially relevant process design in a flowsheet simulator. Therefore, BASF's inhouse flowsheet simulator Chemasim has been chosen, which is now featuring the option to hedge against uncertainties by calculating adjustable robust Pareto sets. To this end, additionally to setting up the MCO problem, the user only distinguishes between HNV and WSV and defines the lower and upper limits of uncertain parameters. The resulting Pareto sets can then be explored interactively [1]. The authors expect that such a transparent handling of uncertainties in process design can contribute to savings in raw materials and energy.

**Acknowledgment**

The authors gratefully acknowledge useful discussions with Roger Böttcher, Kerstin Schneider and Tobias Seidel.

**Entry for the Table of Contents (Graphical abstract)**

**Title:** Pareto-optimal treatment of uncertainties in model-based process design and operation

Jan Schwientek, Katrin Teichert, Jan Schröder, Johannes Höller, Norbert Asprion, Pascal Schäfer2, Martin Wlotzka, and Michael Bortz*

**Type of Article:** Research article

We describe how to deal with here-and-now (i.e. process design related) and wait-and-see (i.e. process operation related) decisions in a multicriteria framework. This approach exploits adjustability for the wait-and-see variables while at the same time respecting optimality guarantees on process KPIs.

Figure TOC

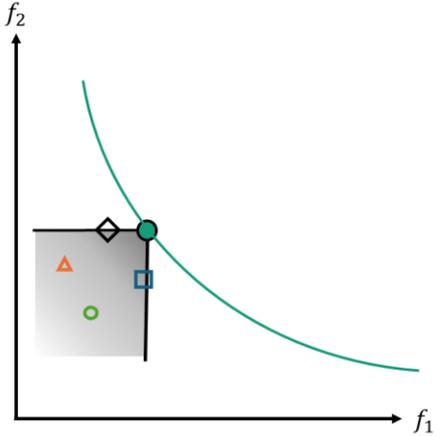

12